# On the use of the variable change *w=exp(u)* to establish novel integral representations and to analyze the Riemann $\varsigma(s,a)$-function, incomplete gamma-function, hypergeometric F-function, confluent hypergeometric $\Phi$-function and beta function

Sergey K. Sekatskii

The variable change $w = e^u$ is applied to establish novel integral representations of the incomplete gamma-function, hypergeometric F-function, confluent hypergeometric $\Phi$-function and beta-function, and to analyze these functions as well as the Riemann $\varsigma(s,a)$-function. In particular, using these representations we give a "pedagogically instructive" proof of the well known approximate functional relation for the Riemann $\varsigma$-function and derive Hurwitz representation of the $\varsigma(s,a)$-function.

# 1. Introduction

In the middle of the XIXth century Schläfli used the variable change $w = e^u$ for the study of the Bessel integral $J_\nu(z) = \dfrac{1}{2\pi i}\int\limits_{-\infty}^{(0+)} w^{-\nu-1}\exp\left(\dfrac{1}{2}z(w-\dfrac{1}{w})\right)dw$ which enabled him to establish some important properties of Bessel functions; see e.g. [1]. Integrals resembling that given above are used to express the Riemann $\varsigma$-function: it is well known that for Re$s$>0:

$$(1-2^{1-s})\Gamma(s)\varsigma(s) = \int\limits_0^\infty \frac{w^{s-1}}{e^w+1}dw \qquad (1)$$

(see e.g. [2] for proof and standard definitions), and similarly for Re$s$>1

$$\Gamma(s)\varsigma(s) = \int\limits_0^\infty \frac{w^{s-1}}{e^w-1}dw \qquad (2)$$

and more general

$$\Gamma(s)\varsigma(s,a) = \int\limits_0^\infty \frac{w^{s-1}e^{aw}}{e^w-1}dw \qquad (3).$$

Similar integral representation is also a starting point of the definition of the incomplete gamma-function [3]

$$\gamma(s,x) = \int\limits_0^x e^{-w}w^{s-1}dw \qquad (4)$$

and actually they appear in many other instances. Hence of course it looks interesting to try the same variable change for analysis of these integrals.

Recently, this approach has been attempted for the Riemann $\varsigma$-function by Bergstrom [4] (the present author is unaware of other researches of this type). In the present Note we analyze some results of its application to different analytical functions.

## 2. Application of the variable change for the Riemann $\varsigma$-function

A variable change $w = e^u$ converts the equality (1) into

$$(1-2^{1-s})\Gamma(s)\varsigma(s) = \int\limits_{-\infty}^\infty \frac{e^{su}}{e^{e^u}+1}du \qquad (5)$$

The integrand $\dfrac{e^{su}}{e^{e^u}+1}$ is a single-valued function on the whole complex plane and thus Cauchy theorem is applicable. This function has simple poles at $u_{n,m} = \ln(\pi(2n+1)) + i\pi(1/2+m)$ where $n$ is integer or equal to zero and $m = 0,\pm 1,\pm 2,...$ The corresponding residues are equal to $-e^{(s-1)u_{n,m}} = i(-1)^m(2n+1)^{s-1}\pi^{s-1}\exp(i\pi s/2)\exp(i\pi ms)$ (because $\frac{d}{du} e^{e^u} |_{u=u_{s,m}} = e^{e^u} \cdot e^u |_{u=u_{n,m}} = -e^{u_{n,m}}$). Let us apply Cauchy theorem to the integral $\int_C \frac{e^{su}}{e^u + 1} du$ taken around the rectangular contour formed by the points $-X$, $\ln(2\pi N)$, $\ln(2\pi N) + 2\pi i$, $-X + 2\pi i$ where $N$ is a positive integer and real $X>0$. Inside this contour we have poles with $m=0, 1$ and $n<2N+1$, whose sum is equal to $S_N = 2\sin(\pi s/2) e^{i\pi s} \pi^{s-1} \sum_{n=0}^{N} (2n+1)^{s-1}$.

Important observation of Bergstrom is that the integral over $(-X + 2\pi i, \ln(2\pi N) + 2\pi i)$ is equal to $e^{2is\pi}$ times the same integral over $(-X, \ln(2\pi N))$ and hence

$$(1 - e^{2is\pi})(1 - 2^{1-s})\Gamma(s)\varsigma(s) = 4\pi i \sin(\pi s/2) e^{i\pi s} \pi^{s-1} \sum_{n=0}^{N} (2n+1)^{s-1} - I_2 + I'_2 \qquad (6).$$

where $I_2$, $I'_2$ are "vertical" integrals taken over $(2\pi N, 2\pi N + 2\pi i)$ and $(-X + 2\pi i, -X + 2\pi i)$ respectively. Writing (6) we have neglected the truncation error (integration over $(-X, \ln(2\pi N))$ instead of $(-X, \infty)$) which is $O(e^{-2\pi N})$. For Re$s>0$ integral $I'_2$ is also trivially negligible when $X \to \infty$. Using $\varsigma(s) = 2^s \pi^{s-1} \sin(\pi s/2)\Gamma(1-s)\varsigma(1-s)$, $\Gamma(s)\Gamma(1-s) = \frac{\pi}{\sin \pi s}$ and $\frac{1 - e^{2i\pi s}}{\sin \pi s} = -2ie^{i\pi s}$ one can thus rewrite

$$-i(1 - 2^{1-s}) e^{i\pi s} \varsigma(1-s) 2^{s+1} \pi^s \sin(\pi s/2) = 4\pi i \sin(\pi s/2) e^{i\pi s} \pi^{s-1} \sum_{n=0}^{N-1} (2n+1)^{s-1} - iI_2 \quad (7)$$

where

$$I_2 = \int_0^{2\pi} \frac{e^{s(\ln(2\pi N)+iy)}}{e^{e^{\ln(2\pi N)+iy}} + 1} dy = 2^s \pi^s N^s \int_0^{2\pi} \frac{e^{isy}}{e^{2\pi N(\cos y + i\sin y)} + 1} dy \qquad (8).$$

We are interested in the range $0 < \sigma < 1$ (as usual $\sigma = \text{Re}\, s$) hence (7) can be rewritten as

$$(1 - 2^{s-1})\varsigma(1-s) = \sum_{n=0}^{N-1} (2n+1)^{s-1} - \tilde{I}_2 \qquad (7a)$$

where

$$\tilde{I}_2 = \frac{2^{s-2} N^s}{\sin(\pi s/2) e^{i\pi s}} \int_0^{2\pi} \frac{e^{isy}}{e^{2\pi N(\cos y + i\sin y)} + 1} dy \qquad (8a).$$

As we mentioned before, equations (7a)-(8a) are not exact but contain an exponentially small error term due to the integration over $(-X, \ln(2\pi N))$ instead of $(-X, \infty)$ (the limitation of integer $N$ is not principal and can easily be removed). Hence, strictly speaking, they are not integral representations of the Riemann $\varsigma$-function although we might consider them as "approximate" integral representations. Note also that they are proven here for Re$s$>1 but by analytical continuation similar relation holds for all $s$ except, of course, $s=1$. For Re$s$<0 integral (8a) vanishes in the limit of large $N$ thus giving the well known relation $(1-2^{s-1})\varsigma(1-s) = \sum_{n=0}^{\infty}(2n+1)^{s-1}$.

For completeness, we finish this Section rewriting equations (7a)-(8a) using the formal change $s \to 1-s$:

$$(1-2^{-s})\varsigma(s) = \sum_{n=0}^{N-1}(2n+1)^{-s} - \tilde{I}_2 \tag{7b}$$

where

$$\tilde{I}_2 = -\frac{2^{-s-1}N^{1-s}e^{i\pi s}}{\cos(s/2)}\int_0^{2\pi}\frac{e^{i(1-s)y}}{e^{2\pi N(\cos y + i \sin y)}+1}dy \tag{8b}.$$

*2.1 Proof of the approximate functional equation for the Riemann $\varsigma$-function*

The natural first approximation to the value of integral (8a) is simply to take the integrand equal to zero for $(0, \pi/2)$ and $(3\pi/2, 2\pi)$, where $\cos y > 0$, and equal to $e^{isy}$ for $(\pi/2, 3\pi/2)$, where $\cos y < 0$. In such a manner we obtain that

$$\tilde{I}_2 \cong \frac{2^{s-2}N^s}{se^{i\pi s}\sin(\pi s/2)}(e^{3i\pi s/2} - e^{i\pi s/2}) = \frac{i2^{s-1}N^s}{s} \tag{9}$$

and hence from (4) we can infer, ignoring for a moment the precision of the expression, that $(1-2^{s-1})\varsigma(1-s) \cong -\frac{2^{s-1}N^s}{s} + \sum_{n=0}^{N-1}(2n+1)^{s-1}$. It is convenient to make a formal change $s \to 1-s$ thus writing $(1-2^{-s})\varsigma(s) \cong -\frac{2^{-s}N^{1-s}}{1-s} + \sum_{n=0}^{N-1}(2n+1)^{-s}$. This is the same as

$$(2^s - 1)\varsigma(s) \cong -\frac{N^{1-s}}{1-s} + \sum_{n=0}^{N-1}(n+1/2)^{-s} \tag{10}.$$

By construction, this is limited with Re$s$<1, but again this relation definitely holds also for Re$s$>1 giving in the limit $N \to \infty$ the well known representation [2] $(2^s - 1)\varsigma(s) = \sum_{n=0}^{\infty}(n+1/2)^{-s}$.

Of course, (10) is the approximate functional equation (applied for $(2^s - 1)\varsigma(s)$ rather than $\varsigma(s)$ but this is not principal) which precision (error term) is long established and is given by the well known theorem (P.77 [2])

THEOREM1. We have

$$(2^s - 1)\varsigma(s) = -\frac{x^{1-s}}{1-s} + \sum_{n \leq x}(n - 1/2)^{-s} + O(1/x^\sigma) \qquad (10)$$

uniformly for $\sigma \geq \sigma_0 > 0$, $|t| < 2\pi x/C$, where $C$ is a given constant greater than 1.

This Theorem is not too difficult to prove by complex variable method as it is done in p. 81 of Ref. [2]; however, the method adopted there requires rather "artificial" introduction of the integral $\int_{x-i\infty}^{x+i\infty} z^{-s} \cot \pi z \, dz$ to obtain, applying the residue theorem, the summing $\sum_{n>x} n^{-s}$ (correspondingly, the integral $\int_{x-i\infty}^{x+i\infty} z^{-s} \tan \pi z \, dz$ and summing $\sum_{n>x}(n-1/2)^{-s}$ in our case). For our approach, this same sum appears quite naturally and hence it would be instructive to prove this theorem based on the proposed variable change, at least by pedagogical reasons. This is done below. For simplicity we soften the conditions of the Theorem1 requiring $|t| < 4x/C$ instead of $|t| < 2\pi x/C$.

As follows from above, for the proof we should rigorously estimate the difference between the actual value of the integral $\int_0^{2\pi} \frac{e^{isy}}{e^{2\pi N(\cos y + i \sin y)} + 1} dy$ and the value $\frac{1}{is}(e^{i3\pi s/2} - e^{i\pi s/2})$ obtained when the integrand is approximated as discussed above. For this, let us divide the integration range as $\int_0^{\pi/2} + \int_{\pi/2}^{\pi} + \int_{\pi}^{3\pi/2} + \int_{3\pi}^{2\pi}$. We will denote these integrals respectively $S_1, S_2, S_3, S_4$, and analyze them one by one.

All these integrals are similar hence we will analyze only $S_4$. (Writing for the integrals $S_{2,3}$ the difference between the integrand at hand and that used for an approximate calculation $\frac{e^{isy}}{1+e^{2\pi N(\cos y + i \sin y)}} - e^{isy} = -\frac{e^{isy}}{1+e^{-2\pi N(\cos y + i \sin y)}}$ we easily see its similarity with $S_{1,4}$). For $S_4$ we make a variable change $y = 3\pi/2 + \tilde{y}$ and obtain $S_4 = e^{3\pi i s/2} \int_0^{\pi/2} \frac{e^{is\tilde{y}}}{e^{2\pi N(\sin \tilde{y} - i \cos \tilde{y})} + 1} d\tilde{y}$.

Let us denote $e^{2\pi N \sin \tilde{y}} = a$ and explicitly write the denominator of the integrand as $e^{2\pi N \sin \tilde{y}}(\cos(2\pi N \cos \tilde{y}) - i \sin(2\pi N \cos \tilde{y})) + 1$. The square of its module is $a^2 + 1 + 2a \cos(2\pi N \cos \tilde{y})$, and now we will demonstrate that always for $\tilde{y} \in [0, \pi/2]$

$$a^2 + 1 + 2a \cos(2\pi N \cos \tilde{y}) \geq a^2/4 \qquad (11).$$

This is easy to see if $a \geq 2$ because then we have $a^2 + 1 + 2a \cos(2\pi N \cos x) \geq a^2 + 1 - 2a = (a-1)^2 \geq a^2/4$, and also (11) is trivially fulfilled if $a \geq 0$ (which is always true) and $\cos(2\pi N \cos \tilde{y}) \geq 0$. When $\cos(2\pi N \cos \tilde{y})$ may be negative? For this one needs $\cos \tilde{y} < 1 - 1/(4N)$ and thus $\sin \tilde{y} = \sqrt{1 - \cos^2 \tilde{y}} \geq 1/(\sqrt{2N})$. But if so, $a = e^{2\pi N \sin \tilde{y}} > e^{\pi \sqrt{2N}} \gg 2$ which finishes the demonstration of (11).

Thus we have (here and below $T = |t|$):

$$\left| \int_0^{\pi/2} \frac{e^{is\tilde{y}}}{e^{2\pi N(\sin \tilde{y} - i \cos \tilde{y})} + 1} d\tilde{y} \right| < 2 \int_0^{\pi/2} \frac{e^{T\tilde{y}}}{e^{2\pi N \sin \tilde{y}}} d\tilde{y} < 2 \int_0^{\pi/2} \frac{e^{T\tilde{y}}}{e^{4N\tilde{y}}} d\tilde{y} =$$

$$\left| \frac{2}{T - 4N} e^{(T-4N)\tilde{y}} \Big|_0^{\pi/2} \right| = \frac{2}{|T - 4N|} |(e^{(T-4N)\pi/2} - 1)| \qquad (12).$$

Here we used $\sin \tilde{y} \geq 2\tilde{y}/\pi$ valid for $\tilde{y} \in [0, \pi/2]$. For negative $T-4N$ (12) is clearly of the order of $O(1/N)$. The factor $e^{3i\pi s/2}$ in front of the integral $S_4$ is cancelled by the denominator $e^{i\pi s} \sin(\pi s/2)$ in front of $\tilde{I}_2$ (see eq. 8a) and thus we have that the error arising from the approximation of all integrals $S_{1,2,3,4}$ is of the order of $O(1/N)$. This finishes the proof: change from $2\pi N$ to arbitrary $x$ is standard for this Theorem, cf. [2], and does not pose any problem.

*2.2. Hurwitz representation of the $\varsigma(s,a)$ function*

Similar consideration can be given for integrals (2) and (3) leading to known integral representations of the functions involved. As an example, let us see how Hurwitz representation of the $\varsigma(s,a)$ [5] can be obtained by the method under discussion. The variable change at question transforms (3) into $\Gamma(s)\varsigma(s,a) = \int_{-\infty}^{\infty} \frac{e^{us} e^{ae^u}}{e^{e^u} - 1} du$. The simple poles are located at $u_{n,m} = \ln(2\pi n) + i\pi(1/2 + m)$ with $n$ integer and $m = 0, \pm 1, \pm 2,...$ Similarly as above, the residues are equal to

$e^{(s-1)u_{n,m}}e^{(-1)^m 2\pi nai} = i(-1)^{m+1}(2n)^{s-1}\pi^{s-1}\exp(i\pi s/2)\exp(i\pi ns)\exp((-1)^m 2\pi nai)$. Application of the same contour as above but with the vertices at $-X$, $\ln(2\pi N+\pi)$, $\ln(2\pi N+\pi)+2\pi i$, $-X+2\pi i$ gives

$$(1-e^{2is\pi})\Gamma(s)\varsigma(s,a) =$$
$$-2\pi i(2\pi)^{s-1}\exp(i\pi s/2)\sum_{n=1}^{N}n^{s-1}(-\exp(2\pi nai)+\exp(-2\pi nai)\exp(i\pi s)) - I_2 + I'_2.$$

(To avoid misunderstanding, under the sum sign we explicitly write the sum of residues for the pair $m=0,1$ for any $n$). Using again $\Gamma(s) = \dfrac{\pi}{\sin\pi s\,\Gamma(1-s)}$ and $\dfrac{1-e^{2i\pi s}}{\sin\pi s} = -2ie^{i\pi s}$ after trivial algebra one can thus rewrite

$$\varsigma(s,a) = 2\Gamma(1-s)(2\pi)^{s-1}\sin(\pi s/2)\sum_{n=1}^{\infty}n^{s-1}\sin(2\pi na+\pi s/2) \qquad (13).$$

where we took the limit $N\to\infty$ hence for Re$s<0$ $\lim_{N\to\infty}I_2 = 0$.

Relation (13) is, of course, Hurwitz formula – again, it is obtained here in a direct and straightforward way without introduction of any "artificial" contour, cf. [5]. In the same way an approximate functional equation for $\varsigma(s,a)$ also can be easily obtained.

### 3. Integral representations for the incomplete gamma-function

As was already mentioned, incomplete gamma-function is defined, whenever appropriate, as $\gamma(s,x) = \int_0^x e^{-w}w^{s-1}dw$ [3] and by analytical continuation otherwise. We make the same variable change $w = e^u$ and note that the integrand contains no poles. Now we take the "more general" rectangular contour formed by the points $-X$, $\ln x$, $\ln x+2\pi ni$, $-X+2\pi ni$ where $n = \pm 1, \pm 2,\ldots$ and for a moment $x$ is positive. Using the same notation as above, we obtain $(1-e^{2is\pi n})\gamma(s,x) = -I_2 + I'_2$ where $I_2 = i\int_0^{2\pi n}\dfrac{e^{s(\ln x+iy)}}{e^{e^{\ln x+iy}}}dy = ix^s\int_0^{2\pi n}\dfrac{e^{isy}}{e^{x(\cos y+i\sin y)}}dy$.

Thus we have a number of representations

$$(1-e^{2is\pi n})\gamma(s,x) = -ix^s\int_0^{2\pi n}e^{isy}e^{-x(\cos y+i\sin y)}dy =$$
$$-ix^s\int_0^{2\pi n}e^{-x\cos y}(\cos(sy-x\sin y)+i\sin(sy-x\sin y))dy \qquad (14).$$

which represent incomplete gamma-function if the values of *sn* are non-integer and which, by analytical continuation, are valid for arbitrary complex values of *x*.

Of course, these representations are closely connected with the known representation [3]

$$\gamma(s,x) = \frac{x^s}{\sin \pi s} \int_0^\pi e^{x\cos\vartheta} \cos(s\vartheta + x\sin\vartheta) d\vartheta \qquad (15)$$

which is usually proven in a rather indirect way using Hankel's generalization of $\Gamma(z)$ as a contour integral, see e.g. p. 244 of Ref. [6]. At the same time, the representation (15) can be easily and directly obtained by our method as follows. Let us introduce an auxiliary function $\tilde{\gamma}(s,x) = \int_0^x e^w w^{s-1} dw$, make the same variable change $w = e^u$ and consider a rectangular integration contour formed by the vertices $-X - \pi i$, $\ln x - \pi i$, $\ln x + \pi i$, $-X + \pi i$ with real $X \to \infty$. It is easy to see that along the line $(-\infty - \pi i, \ln x - \pi i)$ our integral $\int e^{e^z} e^{sz} dz$ is equal to $\int e^{-e^u} e^{su - i\pi s} du = e^{-i\pi s} \gamma(s, x)$ while along the line $(\ln x + \pi i, -\infty + \pi i)$ it is equal to $-e^{i\pi s} \gamma(s, x)$. Thus, repeating all what have been done above (note that the integral $I'_2$ is again negligible for Re*s*>0), we have:

$$(e^{-i\pi s} - e^{is\pi})\gamma(s,x) = -ix^s \int_{-\pi}^\pi e^{isy} e^{x(\cos y + i\sin y)} dy = -2ix^s \int_0^\pi e^{x\cos y}(\cos(sy) + x\sin y) dy.$$ This is exactly the representation (15), and such a representation also can be generalized by taking the rectangular contour formed by the points $-X - n\pi i$, $\ln x - n\pi i$, $\ln x + n\pi i$, $-X + n\pi i$ where *n* is an integer. This, if *sn* is not an integer, gives

$$\gamma(s,x) = \frac{x^s}{\sin \pi n s} \int_0^{\pi n} e^{x\cos\vartheta} \cos(s\vartheta + x\sin\vartheta) d\vartheta \qquad (16).$$

REMARK. When analyzing the Riemann $\varsigma(s,a)$-function, a more general rectangular contour formed by the points $-X$, $\ln(2\pi N)$, $\ln(2\pi N) + 2\pi n i$, $-X + 2\pi n i$, where $n = \pm 1, \pm 2,...$ also can be taken. Quite similarly, one can also repeat the same hint which was used during the analysis of the incomplete gamma-function and take the contour with the vertices at $-X - \pi n i$, $\ln(2\pi N) - \pi n i$, $\ln(2\pi N) + \pi n i$, $-X + \pi n i$, where $n = 1, 2,...$ Unfortunately, this does not lead to any new interesting relations.

## 4. Integral representations for the confluent hypergeometric Φ-function, beta function and hypergeometric F-function

Our next example deals with the integral representation of confluent hypergeometric Φ-function. We have for positive $x$ and $\mathrm{Re}\,c > \mathrm{Re}\,a > 0$ (Ref. [5] P. 255):

$$\Phi(a,c;x) = \frac{\Gamma(c)}{\Gamma(a)\Gamma(c-a)} \int_0^1 e^{xw} w^{a-1} (1-w)^{c-a-1} dw \qquad (17)$$

Our variable change converts this into

$$\Phi(a,c;x) = \frac{\Gamma(c)}{\Gamma(a)\Gamma(c-a)} \int_{-\infty}^0 e^{xe^u} e^{ua} (1-e^u)^{c-a-1} du \qquad (18).$$

Now we take similar rectangular contour as above with the vertices at $-X$, $0$, $2\pi i$, $-X + 2\pi i$ where real $X \to \infty$ and we indent the corners $0$, $2\pi i$ of the contour. Clearly, for $\mathrm{Re}\,c > \mathrm{Re}\,a$ the contributions of integrals over the infinitesimally small quarter-circles related with the indentation vanish, and for $\mathrm{Re}\,a > 0$ the vertical integral $I_2'$ also vanishes. Inside the contour we have no poles and our integrand is a single-valued function, hence we obtain the following new integral representation of the confluent hypergeometrical function:

$$(1 - e^{2\pi i a})\Phi(a,c;x) = -\frac{i\Gamma(c)}{\Gamma(a)\Gamma(c-a)} \int_0^{2\pi} e^{xe^{iy}} e^{aiy} (1-e^{iy})^{c-a-1} dy.$$ This can be rewritten as

$$(1 - e^{2\pi i a})\Phi(a,c;x) =$$
$$-\frac{i\Gamma(c)}{\Gamma(a)\Gamma(c-a)} \int_0^{2\pi} e^{x\cos y} (\cos(x\sin y + ay) + i\sin(x\sin y + ay))(1 - \cos y - i\sin y)^{c-a-1} dy$$

(19).

Similarly as above, taking more general contour $-X$, $0$, $2n\pi i$, $-X + 2n\pi i$ with $n = \pm 1, \pm 2, ...$ we obtain also the representation

$$(1 - e^{2n\pi i a})\Phi(a,c;x) =$$
$$-\frac{i\Gamma(c)}{\Gamma(a)\Gamma(c-a)} \int_0^{2\pi n} e^{x\cos y} (\cos(x\sin y + ay) + i\sin(x\sin y + ay))(1 - \cos y - i\sin y)^{c-a-1} dy$$

(20).

By analytical continuation (19) and (20) are valid for all complex $x$, $a$, $c$ provided $\mathrm{Re}\,c > \mathrm{Re}\,a$ and when $\pm na$ are not integers; of course, the poles of the gamma-function also should be avoided as it takes place also for (18).

The same hint that we used in the previous section can be used to obtain more symmetric integral representation. For this we introduce an auxiliary integral $\int_0^1 e^{-xw} w^{a-1}(1+w)^{c-a-1} dw$, use the same variable change and the contour with the vertices at $-X-\pi n i, -\pi n i, \pi n i, -X+\pi n i$ with real $X \to \infty$.

Our variable change transform the integral into $\int_{-\infty}^0 e^{-xe^u} e^{ua}(1+e^u)^{c-a-1} du$, and it is easy to see that along the line $(-X-\pi n i, -\pi n i)$ it is equal to $e^{-i\pi n a} \int_{-\infty}^0 e^{xe^u} e^{ua}(1-e^u)^{c-a-1} du$ while along the line $(\pi n i, -X+\pi n i)$ it is equal to $-e^{i\pi n a} \int_{-\infty}^0 e^{xe^u} e^{ua}(1-e^u)^{c-a-1} du$ hence

$$2\sin(\pi a n)\Phi(a,c;x) =$$

$$\frac{i\Gamma(c)}{\Gamma(a)\Gamma(c-a)} \int_{-\pi n}^{\pi n} e^{-x\cos y}(\cos(x\sin y - ay) - i\sin(x\sin y - ay))(1+\cos y + i\sin y)^{c-a-1} dy$$

(21).

Here $\pm na$ are not integers and again the poles of the gamma-function should be avoided.

Quite similarly the integral equality defining the beta function (Re$x$>0, Re$y$>0) [Ref. 5 P. 9]:

$$B(x,y) = \int_0^1 w^{x-1}(1-w)^{y-1} dw \qquad (22)$$

can be treated. This leads to the representation

$$(1-e^{2ni\pi x})B(x,y) = -i\int_0^{2\pi n}(\cos(x\vartheta) + i\sin(x\vartheta))(1-\cos\vartheta - i\sin\vartheta)^{y-1} d\vartheta \qquad (23).$$

By analytical continuation, (23) is valid if Re$y$>0 and $\pm nx$ is not an integer. Note the difference between (23) and the integral appearing after the simple variable change $w = \sin\theta$ (or similar) in (22):

$$B(x,y) = \int_0^{\pi/2} \sin^{x-1}\vartheta(1-\sin\vartheta)^{y-1} \cos\vartheta d\vartheta \qquad (24).$$

Again, the same hint that was used for an incomplete gamma-function to obtain more symmetric integral representation can be repeated by introducing an auxiliary integral $\int_0^1 w^{x-1}(1+w)^{y-1} dw$ and considering the contour formed by the vertices $-X-\pi n i, -\pi n i, \pi n i, -X+\pi n i$ with real $X \to \infty$:

$$2\sin(\pi n x)B(x,y) = \int_{-n\pi}^{n\pi}(\cos(x\vartheta) + i\sin(x\vartheta))(1+\cos\vartheta + i\sin\vartheta)^{y-1} d\vartheta \qquad (25)$$

where $nx$ is not an integer.

The same approach can be used to obtain an integral representation for the hypergeometric function $F(a,b;c;z)$ starting from the known representation $F(a,b;c;z) = \dfrac{\Gamma(c)}{\Gamma(b)\Gamma(c-b)} \int_0^1 w^{b-1}(1-w)^{c-b-1}(1-wz)^{-a} dw$ valid for $\operatorname{Re} c > \operatorname{Re} b > 0$, $|\arg(1-z)| < \pi$; [5] P. 114. For these conditions our standard rectangular contour with indented corners which we used to study the $\Phi(a,c;x)$ function again contains no poles and we have

$(1 - e^{2n\pi i b}) F(a,b;c;z) =$
$-\dfrac{i\Gamma(c)}{\Gamma(a)\Gamma(c-a)} \int_0^{2\pi i} (\cos(by) + i\sin(by))(1 - \cos y - i\sin y)^{c-b-1}(1 - z\cos y - iz\sin y)^{-a} dy$

(26).

Here $nb$ is not an integer and the poles of gamma function should be avoided. Similarly as above, more symmetrical representation can be obtained introducing an auxiliary integral
$\int_0^1 w^{b-1}(1+w)^{c-b-1}(1+wz)^{-a} dw$ and considering the contour formed by the vertices $-X - \pi n i$, $-\pi n i$, $\pi n i$, $-X + \pi n i$ with real $X \to \infty$:

$2\sin(\pi n b) F(a,b;c;z) =$
$\dfrac{\Gamma(c)}{\Gamma(a)\Gamma(c-a)} \int_{-n\pi}^{n\pi} (\cos(by) + i\sin(by))(1 + \cos y + i\sin y)^{c-b-1}(1 + z\cos y + iz\sin y)^{-a} dy$ (27).

## 5. Conclusion

Thus we have shown that the variable change $w = e^u$ is useful for an analysis of integrals representing the Riemann $\varsigma$–function, incomplete gamma-function, confluent hypergeometric $\Phi$-function and beta function. There are no doubts that other interesting and yet unexplored applications of this variable change can be found.

And now let us make the final remark. Might it be reasonable to propose the introduction of the "incomplete Riemann $\tilde{\varsigma}$-function" which may be defined for positive $x$ and $\operatorname{Re} s > 1$ as $\Gamma(s)\tilde{\varsigma}(s,x) = \int_0^x \dfrac{w^{s-1}}{e^w - 1} dw$ or similarly and for other parameters by an appropriate analytical continuation? This function deepens the analogy between the Riemann function and gamma function and might be useful for some physical applications. For example, it is well known that the spectral density of thermal body radiation energy is expressed in 3D via the Riemann function $\Gamma(4)\varsigma(4) = \int_0^\infty \frac{w^3}{e^w - 1} dw$ (and via other related functions in other number of dimensions and may be in fractals also (?)). The highest possible frequency of the system can be limited by different physical reasons (e.g. by Planck length and time scales if nothing else) which apparently makes the introduction of the corresponding "incomplete" function quite reasonable: in 3D the real spectral density of the blackbody radiation is $\Gamma(4)\tilde{\varsigma}(4, x) = \int_0^x \frac{w^3}{e^w - 1} dw$, and so on.

REFERENCES
[1] N. Watson, A treatise on the theory of Bessel functions, Cambridge, Cambridge Univ. Press, 1958, P. 176.
[2] E. C. Titchmarsh and E. R. Heath-Brown, The theory of the Riemann Zeta-function, Oxford, Clarendon Press, 1988.
[3] A. Erdelyi, Higher transcendental functions, Vol. II, McGrow and Hill, N-Y., 1953, P. 137.
[4] A. Bergstrom, Proof of Riemann's zeta-hypothesis, arXiv:0809.5120, 2008.
[5] A. Erdelyi, Higher transcendental functions, Vol. I, McGrow and Hill, N-Y., 1953, P. 26.
[6] E. T. Whittaker and G. N. Watson, A course of modern analysis, Cambridge, Cambridge Univ. Press, 1965.

S. K. Sekatskii, Laboratoire de Physique de la Matière Vivante, IPSB, BSP 408, Ecole Polytechnique Fédérale de Lausanne, CH 1015 Lausanne-Dorigny, Switzerland. E-mail : serguei.sekatski@epfl.ch